\documentclass[%
english,%
twoside
]{%
amsart%
}

\usepackage{amscd}

\usepackage{amssymb}

\usepackage{a4}

\newtheorem{statement}{}

\newcommand{\ep}{\hspace*{\fill}$ \lozenge $}

\newcommand{\normal}{\vartriangleleft}

\newcommand{\opon}{\ltimes}

\newcommand{\C}{{\mathbb C}}
\newcommand{\Z}{{\mathbb Z}}

\newcommand{\be}{\begin{equation}}
\newcommand{\ee}{\end{equation}}
\newcommand{\bea}{\begin{eqnarray}}
\newcommand{\eea}{\end{eqnarray}}
\newcommand{\beann}{\begin{eqnarray*}}
\newcommand{\eeann}{\end{eqnarray*}}
\newcommand{\bs}{\begin{statement}}
\newcommand{\es}{\end{statement}}

\DeclareMathOperator{\Ad}{Ad}
\DeclareMathOperator{\Aut}{Aut}
\DeclareMathOperator{\Der}{Der}
\DeclareMathOperator{\End}{End}

\DeclareMathOperator{\Hom}{Hom}
\DeclareMathOperator{\Id}{Id}
\DeclareMathOperator{\Image}{Image}
\DeclareMathOperator{\Ker}{Ker}

\title{%
	Generalised Magnus modules
	over the braid group
      }
\author{%
	Mirko L\"{u}dde
}
\address{%
	Graduiertenkolleg
	``Geometrie und Nichtlineare Analysis''
\\
	Institut f\"ur Reine Mathematik
\\
	Humboldt Universit\"at zu Berlin
\\
	Ziegelstrasse 13a
\\
	D-10099 Berlin, Germany
	}
\thanks{%
	Supported by {\sl Deutsche Forschungsgemeinschaft}.%
\\
	A prior version appeared as preprint 170
	of the ``SFB 288'' in June 1995.
	See ``http://www.math.tu-berlin.de''%
}
\email{%
	luedde@mathematik.hu-berlin.de
}
\date{%
	July 10, 1995
}
\subjclass{%
	Primary		20F36;	
	Secondary	32S40	
}
\keywords{%
	Braid Group, Magnus Representation, Burau Representation%
}
\begin{document}
\begin{abstract}
	W.~Magnus' representations of submonoids
	$ E \leq \mbox{End}(F) $
	of the endomorphisms of a free group
	$ F $ of finite rank are generalised
	by identifying them with the first homology group
	of $ F $ with particular coefficient modules.
	By considering
	a suitable free resolution of the integers
	over the semidirect product of free groups,
	a class of representations of the braid group
	can be obtained on higher homology groups.
	The resolution shows that
	the holonomy representations
	of the braid group and of the Hecke algebra
	constructed topologically by R.~J.~Lawrence
	belong to this class.
\end{abstract}

\maketitle

\tableofcontents

\section{%
	Introduction
	}
These notes were prepared for seminar talks.
They are based on some of the ideas sketched
in the author's thesis \cite{Luedde1992}.
It is shown that the class of Magnus representations
of the braid group,
cf.\cite{Birman},
can be suitably generalised
via group homology
to include the topological representations of
R.~J.~Lawrence \cite{Lawrence1990}.

W.~Magnus, in order to give a derivation for the
well-known Burau representation of the braid group,
described a class of representations
of automorphism groups of free groups,
see e.g.~\cite{Magnus1974}.
J.~Birman devoted a chapter of her book
to this class of
``Magnus representations'',
cf.\cite{Birman}.
Among these representations are
e.g.\,the representation
found by W.~Burau,
which leads to a construction of the
Alexander polynomial via a Markov functional,
as well as those constructed by B.~Gassner.
G.~D.~Mostow implicitly used Magnus representations
of the unpermuting braid group
to investigate the monodromy group
of Euler-Picard integrals
\cite{Mostow1987}.
This is just to mention some of the applications.

The Magnus representation modules
(in the sense originally considered by W.~Magnus)
can be understood as the first homology groups of the free group
with particular coefficient modules.
Using Artin's imbedding
$
B_n \rightarrow \Aut (F_n)
$,
a semidirect product
$
B_n \opon F_n
$
of the braid group $ B_n $ with
the free group $ F_n $ of rank $n$
can be defined.
$F_n$ is imbedded as a normal subgroup
into this product.
By functoriality of $H_*$ the quotient group
$
( B_n \opon F_n ) / F_n \simeq B_n
$
acts onto
$ H_1(F_n,M), $
if $M$ is a module over
$B_n \opon F_n$.
We will show that this construction generalises
the approaches of Magnus and Birman.

These considerations show that Magnus' representations are
the first nontrivial examples of
a general procedure:
one may use any pair $K \normal G$ of groups with
quotient $G / K \simeq B_n$ and any $G$ module $N$
to obtain a $B_n$ module $H_*(K,N)$.

A step in this direction has been undertaken by
R.~J.~Lawrence, cf.\cite{Lawrence1990},
from a different and geometrical point of view.
In an attempt to understand the Jones polynomial
in geometrical terms,
she constructs a vector bundle with a
natural flat connection over a
base space having
the braid group $B_n$ as its fundamental group.
In this way she obtains braid representations from the holonomy
of the connection,
which under certain conditions factor through
the Hecke algebra.
The typical fiber of the bundle is the $m$-th homology
$ H_m(Y; \chi ) $
of the configuration space $Y$
of $m$ distinct points in the $n$-fold
punctured plane with a suitably chosen abelian local
coefficient system
$ \chi \in \Hom (\pi_1(Y,y), \C \backslash \{ 0 \} ) $.
The class of representations obtained in this way e.g.
can be used for the construction of the one-variable
Jones polynomial \cite{Lawrence111990}.
For a short account on this approach,
cf.\cite{Atiyah1989}.
The case $m = 1$ was investigated
with different motivation in
\cite{Mostow1987}.

In the present note, it will be shown that,
since $Y$ is an Eilenberg-MacLane complex of type
$ (\pi _1(Y,y), 1), $
Lawrence's construction
via Eilenberg's theorem on local coefficients
is the topological approach
to the construction of the homology of the group
$ \pi _1(Y,y) $ with coefficient module
$ (\chi , \C ). $
The fundamental group
$ \pi _1(Y,y) $
can be imbedded as a normal subgroup
into a generalised braid group
$ B_{n,m} $
with quotient
$ B_{n,m} / \pi _1(Y,y) \simeq B_n. $
Therefore Lawrence's approach
fits into the algebraic setting
sketched above.
In particular
$ \pi_1(Y,y) \simeq F_n $
for $m = 1$
and we recover Magnus' representations.

We will therefore,
after having reviewed Magnus' original construction,
find a free $\pi _1(Y,y)$ resolution
$ \partial \in \Hom _{\pi _1(Y,y)}(C) $
of the integers.
This resolution is supposed by the recursive structure of
$ \pi _1(Y,y). $
It explicitly allows the computation of the $B_n$
action onto chain complexes
$ N \otimes _{\pi _1(Y,y)} C $
by means of braid ring valued matrices.
A special version of these matrices was found
topologically in \cite{Lawrence1990}.
In their general form they are
the ``braid valued Burau matrices''
derived independently and by different means
in \cite{ConstantinescuLuedde1992}, \cite{Luedde1992}.
The name is due to the fact
that the classical Burau matrices
are images of the aforementioned ones
under a ring homomorphism of the matrix elements.
The braid valued matrices encode
all the braid representations mentioned
above in a simple and unified form.

{\bf Acknowledgements:}
This work has been made possible
by support from the
{\sl Graduiertenkolleg
Geometrie und Nichtlineare Analysis}
at the
{\sl Humboldt Universit\"at zu Berlin}.
It is a pleasure for me to thank
Prof.~Th.~Friedrich and the
{\sl Deutsche Forschungsgemeinschaft}.
I also thank the participants
(in particular A.~Nestke) of our
seminar on Vassiliev invariants and related topics
for their interest and for clarifying discussions
and F.~Neumann for sending me manuscripts prior
to publication.

\section{%
	Magnus representations of submonoids of $\End (F)$
        }
\label{SectionMagnus}
This section presents a description
of the ``classical'' situation
in a more general formulation.
The approaches of
\cite{Birman}
and
\cite{Magnus1974}
will be linked together.

In the following,
$E$ is a general submonoid
of the monoid $\End (F)$ of right endomorphisms
of the free group $F$.
A subscript indicates the finite rank $n$ of $F$,
if necessary.
Since the endomorphisms are acting from the right,
the product $ab$ is the endomorphism $b \circ a$,
$b$ after $a$.
The semidirect product
$ E \opon F $
is the monoid consisting of the set
$ E \times F $
equipped with the product
$
(a,f)(b,g) = (ab, b(f)g),
$
$a,b \in E$, $f,g \in F$
and with unity $(1,1)$.
We will identify $E$ and $F$ with their images
in $E \opon F$ and therefore have
$ af = (a,1)(1,f) = (a,f), $
$ E \opon F = EF $
and
$ afbg = abb(f)g. $
$F$ is normal in $EF$
in the sense that for every $p \in EF$
and $f \in F$ there is a $g \in F$
uniquely determined by $fp = pg$.
The linear extension
of any group (or monoid) homomorphism
$
\phi \in \Hom (G,H)
$
to the group ring $\Z G$
is the map $\phi _*$.
The relative augmentation ideal
$
I_{N,G} = \Ker (\phi _*) \normal \Z G
$
is generated by the subset
$
N - 1 \subset \Z G
$,
where $N$ is the normal subgroup
$
N = \Ker (\phi ) \normal G
$.
We set
$I_{G,G} = I_G$
if we consider full augmentation ideals.

\begin{statement}[W.~Magnus]
Let
$U \normal F$
be a normal subgroup.
Let
$
E = \mbox{St}_{\End (F)}(F / U) \leq \End (F)
$
be the stabilizing monoid
of $F/U$ in $\End (F)$.
Then $U/[U,U]$ is an
$F/U$-$E$ bimodule:
$[f]\in F/U$ maps $[u]\in U/[U,U]$
to
$[fuf^{-1}]$, $a \in E$ maps $[u]$ to
$[a(u)]$,
and the actions of $E$ and $F/U$ commute.
\end{statement}
This statement is taken from \cite{Magnus1974}, p.471.
The assignments defined by $[f]$ and $a$ are mappings:
let $u,p\in U$,  $c\in [U,U]$.
Then
$
fpucp^{-1}f^{-1} = fupp^{-1}u^{-1}pucp^{-1}f^{-1}
		= fuf^{-1}c'
$
for some $c'\in [U,U]$.
Therefore the image of $[u]$ under $[f]$ is uniquely
determined.
For $ a\in E $ we have
$ a(uc) = a(u)a(c) $
and
$ a(c) \in [U,U] $,
such that $a$ too defines a map.
Since
$E$
acts trivially
onto $F/U$,
$a(f) = fp$ for some $p \in U$.
Thus,
$
a(fuf^{-1})
	= fpa(u)p^{-1}f^{-1}
	= fa(u)f^{-1}fp [p^{-1},a(u)^{-1}] p^{-1}f^{-1},
$
where
$
fp[p^{-1},a(u)^{-1}]p^{-1}f^{-1} \in [U,U]
$,
and the operations of
$E$ and $F/U$ commute.
\ep

Using the Fox derivation,
Birman has given another construction
related to the Magnus modules $U / [U,U]$.
For any epimorphism
$\psi \in \Hom (F,H)$
let
$ M_{\psi } $
be the free left ${\bf Z}H$ module of rank $n$
with basis
$ \{ s_1,\ldots ,s_n\} $.
Let
$
\delta _{\psi }  \in \Der (F,M_{\psi })
$
be the total Fox derivation,
cf.~\cite{Birman}, pp.104.
It is a group derivation uniquely defined
as a map by the equations
\begin{eqnarray*}
\delta _{\psi }(f_i) & := & s_i,
\\
\delta _{\psi }(fg) & := &
	\delta _{\psi }(f) + \psi (f)\delta _{\psi }(g),
\end{eqnarray*}
for a free generating system
$ \{ f_1,\ldots ,f_n \} $
of $F$ and $f,g \in F$.
$ \delta _{\psi } $ has an extension
$ \bar{\delta }_{\psi } \in \Der (\Z F, M_\psi ) $
to a ring derivation,
defined as the unique homomorphism of abelian groups
with
$ \bar{\delta }_{\psi }(f) = \delta _{\psi }(f), $
$ f \in F $.
We drop the bar distinguishing both derivations.
Let
$
\partial _{j,\psi} \in \Der (F,\Z H)
$
be the partial Fox derivations with respect to the generators
$f_j$,
$
\delta _{\psi }(f) = \sum _{j} \partial _{j,\psi}(f) s_j
$,
$f\in F$.
When we consider the case $\psi = \Id $
we drop indices referring to $\psi $.

\begin{statement}[J.~Birman]
Let $E \leq \End (F)$
be the maximal monoid satisfying
$
\psi \circ a = \psi
$
for all $a \in E$.
For every $a \in E$
let $\rho (a)$ be the endomorphism of $M_{\psi }$
uniquely defined by
$$
s_{i} \mapsto \delta _{\psi }(a(f_{i})).
$$
Then
$
\rho \in \Hom (E, \End (M_{\psi })).
$
\end{statement}
Statement and proof follow
\cite{Birman}, thm.3.9, p.116.
For $a \in E$ let the ``Fox-Jacobi matrix''
with respect to the generators
$ \{ f_i \} $
be defined by
$
a_{i,j} = \partial _j(a(f_i))
$.
Due to the chain rule for the Fox derivation,
\cite{Birman}, prop.3.3, pp.105,
for $a, b \in E$ holds,
\begin{eqnarray*}
(ab)_{i,j} & = & \partial _j(b(a(f_i)))
\\
           & = & \sum _k  b(\partial _k(a(f_i)))
                  \partial _j(b(f_k))
\\
           & = & \sum _k b(a_{i,k})b_{k,j}.
\end{eqnarray*}
Using
$
\rho(a)(s_i)
= \delta _{\psi }(a(f_i))
= \sum _j \partial _{j,\psi }(a(f_i))s_j
= \sum _j \psi _*(a_{i,j})s_j
$,
we find,
\begin{eqnarray*}
\rho(ab)(s_i) & = & \sum _j \psi _*((ab)_{i,j})s_j
\\
                    & = & \sum _{j,k}\psi _*(b(a_{i,k})b_{k,j})s_j
\\
                    & = & \sum _{j,k}\psi _*(a_{i,k})\psi _*(b_{k,j})s_j
\\
                    & = & \rho(b)(\rho(a)(s_i)).
\end{eqnarray*}
Thus,
$
\rho \in \Hom (E,\End (M_{\psi }))
$.
\ep

In order to understand
the relationship between
the Magnus modules $U / [U,U]$
and the Birman modules $M_{\psi }$,
let
$
\mu _{\psi } \in \Hom (F, H \opon M_{\psi })
$
be the
representation of $F$
that is associated to the derivation
$
\delta _{\psi }  \in \Der (F,M_{\psi })
$:
the semidirect product of
$H$ and $M_{\psi }$ is defined by
$
(g,p)(h,q) = (gh, p + g q)
$
for
$g, h \in H$
and
$p, q \in M_{\psi }$
and $\mu _{\psi }$ maps
$
f \mapsto (\psi(f), \delta _{\psi }(f)).
$
Since $\delta _{\psi }$ is a group derivation,
$\mu _{\psi }$ is a homomorphism
(called the ``Magnus $\psi $ representation'' of $F$
in \cite{Birman}).
It is characterised by
\begin{statement}[Blanchfield, Birman]
\begin{eqnarray*}
\Ker (\mu _\psi) & = &
 [\Ker (\psi ), \Ker (\psi )],
\\
\Image (\mu _\psi) & = &
\{ (h, \sum h_j s_j); h \in H, h_j \in \Z H,
h - 1 =  \sum h_j (\psi (f_j) - 1) \} .
\end{eqnarray*}
\end{statement}
By theorems of Blanchfield,
cf.~\cite{Birman}, thm.3.5, pp.107,
and Magnus,
cf.~\cite{Birman}, thm.3.7, pp.111,
the kernel is identified.
The image is determined
by \cite{Birman}, thm.3.6, pp.108.
\ep

Both constructions are related by the following fact.
\begin{statement}
Given either a Magnus module
$ U / [U,U] $
or a Birman module
$ M_{\psi }, $
there is a short exact sequence of
$\Z H$-$E$ bimodules,
$$
0 \rightarrow U / [U,U] \rightarrow M_{\psi } \rightarrow Q \rightarrow 0,
$$
with
$ H \simeq F / U, $
$ \psi \in \Hom (F,H) $
such that
$ \Ker (\psi ) = U, $
and
$ E = \mbox{St}_{\End (F)}(F / U). $
\end{statement}
Let $U \normal F$ be given,
$
E = \mbox{St}_{\End (F)}(F/U),
$
giving rise to the Magnus module
$ U / [U,U] $
over $F / U$ and $E$.
Let
$
\psi \in \Hom (F, F/U)
$
be the projection homomorphism.
Then
$
\psi (a(f)) = \psi (fp) = \psi(f)
$,
with some $p \in U$,
since $E$ acts as identity onto $F/U$.
So the projection
$\psi $ is $E$-invariant
and it defines a Birman module $M_{\psi }$.
We must find a monomorphism
from the $F/U$-$E$ bimodule $U/[U,U]$ into
$M_{\psi }$.
The Magnus representation of $F$, $\mu _{\psi }$,
induces an isomorphism of the group $F/[U,U]$
onto its image.
The restriction of this monomorphism to $U/[U,U]$
satisfies
$ \mu _{\psi }([u]) = (1, \delta _{\psi }(u)) $,
for
$ u \in U = \Ker (\psi ). $
So
$
\delta _{\psi } : U/[U,U] \rightarrow M_{\psi }
$
has the desired properties:
it is an injective homomorphism of abelian groups
due to Magnus' theorem on
$
\Ker (\mu _{\psi })
$.
It intertwines the operation of $\Z (F / U)$ onto
$U/[U,U]$
with the left multiplication in
$M_{\psi }$,
\begin{eqnarray*}
\delta _{\psi }(f u f^{-1}) & = &
\delta _{\psi }(f)
+ \psi(f) \delta _{ \psi }(u)
- \psi(f u f^{-1}) \delta _{\psi }(f)
\\
& = &
\psi (f) \cdot \delta _{\psi }(u).
\end{eqnarray*}
It also intertwines the right $E$ actions on the two modules,
$
\delta _{\psi }(a([u]))= \rho(a)(\delta _{\psi }([u])).
$
Indeed, from
$
u = f_{i_{1}}^{\epsilon _{1}} \ldots f_{i_{k}}^{\epsilon _{k}}
$,
$\epsilon _{j} \in \{ -1, 1 \} $
follows
$$
\delta _{\psi }(u) =
\sum _{j=1}^{k}
\psi (f_{i_{1}}^{\epsilon _{1}} \ldots
f_{i_{j-1}}^{\epsilon _{j-1}})
\epsilon _{j} \psi(f_{i_{j}})^{(\epsilon _{j}-1)/2}
s_{i_{j}}
$$
and
$$
\delta _{\psi }(a(u)) =
\sum _{j=1}^{k}
	\psi (a(f_{i_{1}}^{\epsilon _{1}}
		\ldots f_{i_{j-1}}^{\epsilon _{j-1}}))
	\epsilon _{j}
	\psi(a(f_{i_{j}}^{(\epsilon _{j}-1)/2}))
	\delta _{\psi}(a(f_{i_j})).
$$
Now the $E$-invariance of $\psi $ can be used
and by comparison of
the expressions the claim is proved.

Conversely, let $\psi $ be given and be $E$-invariant.
$ U = \Ker (\psi ) $
is an $E$-invariant
(since
$ \psi (a(u)) = \psi (u) = 1 $
for $ u \in U $)
normal subgroup of $F$.
$E$ acts as identity onto $F/U$,
since
$ f^{-1}a(f) \in U $,
so $ U / [U, U] $ is a Magnus module
with respect to $E$.
Again the Fox derivation
$ \delta _{\psi } \in \Der (U, M_{\psi }) $
yields a monomorphism
$ U/[U,U] \rightarrow M_{\psi } $.
\ep

We call this sequence the
``Birman-Magnus sequence''.
The image
$ \delta _{\psi }(U) $
of the Magnus module
in the Birman module ``measures'' how far the set
$\{ \psi(f_j) - 1 \} $
is from being linearly independent over
${\bf Z}H$.
\begin{statement}\label{TheoremFoxImage}
$
\Image (\delta _{\psi } : \Ker (\psi ) \rightarrow M_{\psi }) =
\{ \sum _j h_j s_j;
h_j \in {\bf Z}H,
\sum h_j (\psi (f_j) - 1) = 0 \} .
$
\end{statement}
Let
$ u \in \Ker (\psi ) $
and
$ \delta _{\psi }(u) = \sum h_j s_j, $
for some
$ h_j \in \Z H. $
Then
$ u - 1 \in I_F $,
where $I_F$ is the augmentation ideal
of $F$ in the ring $\Z F$,
such that
$ u - 1 = \sum u_j (f_j - 1) $
for suitable (and unique)
$ u_j \in \Z F. $
We apply
$ \delta _{\psi } \in \Der (\Z F, M_{\psi }) $
to find
$ \sum h_j s_j  = \sum \psi _*(u_j) s_j $
and
$ h_j = \psi _* (u_j). $
Then
$ \psi _* (u - 1) = 0 $
implies,
$ \sum h_j (\psi (f_j) - 1) = 0. $
Conversely, let
$ \sum h_j (\psi(f_j) - 1) = 0, $
$ h_j \in \Z H. $
By Lyndon's theorem,
\cite{Birman}, thm.3.6, pp.108,
there is a
$ u \in \Ker (\psi ) $
with
$ \delta _{\psi}(u) = \sum h_j s_j. $
\ep

Now we will construct
a short exact sequence of $EF$-$E$ bimodules
that generalises the Birman-Magnus sequence.
Let
$ \eta \in \Hom (EF, E) $
be the projection,
i.e.~$ \eta (af) = a $
for
$ a \in E $,
$ f \in F $.
The extension
$ \eta _* \in \Hom (\Z (EF), \Z E) $
to the group ring
defines the relative augmentation ideal
$
\Ker (\eta _*) = I_{F,EF}
$
of $F$ in $\Z EF$.
We consider the ideal $I_{F,EF}$
as left $EF$ right $E$ bimodule.
Since it is free as a left $EF$ module,
there is a representation of $E$
in terms of $\Z (EF)$-valued matrices.
In section \ref{SectionBurau} we
will find these matrices for the case
of the braid group.

\begin{statement}
\label{TheoremSequence}
Let $N$ be a right $EF$ module.
Then there is a short exact sequence
of right $E$ modules,
$
0
\rightarrow
H_1(F,N)
\rightarrow
N \otimes _{EF} I_{F,EF}
\rightarrow
N I_{F,EF}
\rightarrow
0.
$
\end{statement}
We have a free resolution
$
0
\rightarrow
P_1
\stackrel{\partial _1}{\rightarrow }
P_0
\stackrel{\epsilon _*}{\rightarrow }
\Z
\rightarrow
0
$
of the integers over the free group $F$,
where
$P_1 = I_F$ is the augmentation ideal,
$P_0 = \Z F$ is the group ring and
$\epsilon _*: \Z F \rightarrow \Z $
is the augmentation homomorphism,
cf.~\cite{Robinson}, thm.11.3.2, p.322.
By the exact sequence
$
0 \rightarrow H_1(F,N)
\rightarrow
N \otimes_F \partial_1(P_1)
\rightarrow
N \otimes_F P_{0}
$,
cf.~\cite{Robinson}, thm.11.2.7, p.319,
we have
$
H_1(F, N) \simeq
\Ker (N \otimes _F I_F \rightarrow N \otimes _F \Z F)
$.
Finally we notice that
$
N \otimes _F I_F \simeq N \otimes _{EF} I_{F,EF}
$
(on the right hand side, the structure as
right $E$ module is more obvious)
and that the maps
$
N \otimes _F I_F \rightarrow N \otimes _F \Z F
$
and
$
N \otimes _{EF} I_{F,EF} \rightarrow N I_{F,EF}
$
have isomorphic kernels,
the later map is onto
and the kernel is preserved by the right $E$ action.
Thus, we obtain the short exact sequence
of right $E$ modules as claimed.
\ep

We will show
that this sequence
specialises to
the Birman-Magnus sequence.
Let
$
\psi \in \Hom (F,H)
$
be an epimorphism with
$
\psi \circ a = \psi
$
for every
$ a \in E. $
This condition allows us
to construct an extension
$
\bar{\psi } \in \Hom (EF,H),
$
by setting
$ \bar{\psi }(af) = \psi (f), $
$ a \in E, $
$ f \in F $.
We regard
$ \Z H $
as a right $EF$ module using
$ \bar{\psi }, $
i.e.~$ af: h \mapsto h \bar{\psi }(af) = h \psi (f). $
\begin{statement}
\label{TheoremIdentify}
\begin{enumerate}
\item
	The Birman module
	$M_{\psi }$
	as left $\Z H$ right $E$ bimodule
	is isomorphic to
	$
	\Z H \otimes _{EF} I_{F,EF}.
	$
\item
	Let
	$U = \Ker (\psi ).$
	The corresponding Magnus module
	$ U / [U,U] $
	is isomorphic to
	$ H_1(F, \Z H). $
\end{enumerate}
\end{statement}
1)
The restriction of $\delta _{\psi }$ to $I_F$ is a homomorphism
of left $F$ modules,
$
\delta _{\psi }(f (f_i - 1)) =
\delta _{\psi }(ff_i) - \delta _{\psi }(f) =
\psi (f) \delta _{\psi }(f_i - 1).
$
So there is a map
$$
\Z H \otimes _{EF} I_{F,EF} \simeq \Z H \otimes _F I_F
\stackrel{1 \otimes \delta _{\psi}}{\longrightarrow }
\Z H \otimes _F M_{\psi } \simeq M_{\psi }.
$$
We claim it is an isomorphism of
$\Z H$-$E$ bimodules.
Since
$
\Z H \otimes _F I_F
$
and
$ \Z H \otimes _F M_{\psi } $
are free left $\Z H$ modules with bases
$ \{ 1 \otimes (f_j - 1) \} $
and
$ \{ 1 \otimes s_j \} , $
respectively,
and
$
(1 \otimes \delta _{\psi })(1 \otimes (f_i - 1)) =
1 \otimes s_j
$,
the map is an isomorphism of left $\Z H$ modules.
We now consider the right $E$ action.
In
$ \Z H \otimes _{EF} I_{F,EF} $
for $ a \in E $ we have
$
1 \otimes (f_j - 1)a = 1 \otimes a (a(f_j) - 1)
$
which is equal to
$
1 \otimes (a(f_j) - 1) \in \Z H \otimes _F I_F
$
due to the condition
$ \bar{\psi }(a) = 1. $
By
$ 1 \otimes \delta _{\psi} $
the element
$ 1 \otimes (a(f_j) - 1) $
is mapped to
$
1 \otimes \delta _{\psi }(a(f_j)),
$
such that indeed we have a homomorphism
of right $E$ modules.

2)
The Birman-Magnus sequence together
with theorem (\ref{TheoremFoxImage})
identifies the Magnus module
$ U / [U,U] $
with the submodule
$
\{ \sum _j h_j s_j;
h_j \in {\bf Z}H,
\sum h_j \psi _*(f_j - 1) = 0 \} \leq M_{\psi }
$
via
$ \delta _{\psi }. $
This submodule is isomorphic to
$
\Ker (\Z H \otimes _F I_F \rightarrow \Z H \otimes _F \Z F)
\simeq H_1(F, \Z H)
$
by the map
$
\sum h_j s_j
\mapsto
\sum h_j \otimes (f_j - 1).
$
\ep

We thus have generalised the Birman-Magnus modules
$ \Z H \otimes _{EF} I_{F,EF} $
insofar, as we allow the condition
$ \psi \circ a = \psi $
to be dropped and any
right $ EF $ module $N$ be used instead of
the ring
$ \Z H, $
which is a particular right $ EF $ module.

There are no
higher nontrivial homology groups
than $ H_1(F, N) $.
But by choosing
$ E $ as
Artin's braid group $ B_n $,
one may use semidirect products
of free groups in order to construct
sequences of modules over $B_n$
from higher homology groups.
The homology module discussed so far is
the first member of this sequence.

\section{%
	Generalities on the braid group
}
{}From now on we specialise our considerations
from $E \leq \End (F)$ to Artin's braid group.
In this section we
collect some well known facts
about it.
Eventually we show that
the topological homology groups appearing as the fibers
in the geometrical approach to braid representations in
\cite{Lawrence1990}
can be understood as group homology.

The braid group
$B_{n}$ on $n$ strings
(over the Euclidean plane)
is the
group generated by the $n-1$-set
$
\{ \tau _{i}; i \in \{ 1, \ldots , n-1 \} \}
$
according to the relations of Artin,
cf.\cite{Birman}, lemma 1.8.2, pp.20:
\begin{eqnarray*}
\tau _{i} \tau _{j} & = &
\tau _{j} \tau _{i}
\mbox{, if } \mid i-j \mid \geq 2,
\\
\tau _{i} \tau _{1+i} \tau _{i} & = &
\tau _{1+i} \tau _{i} \tau _{1+i}.
\end{eqnarray*}
We set
$ \tau _{i,i} := \vartheta _{i,i} := 1 $
and use the abbreviations
($ i,j,k,l \in \{ 1,\ldots ,n \} $; $i < j$)
\begin{eqnarray*}
\tau _{i,j} & := &
	\tau _i \tau _{1+i} \ldots \tau_{j-2} \tau _{j-1},
\\
\tau _{j,i} & := &
	\tau _{j-1} \tau _{j-2} \ldots \tau_{1+i} \tau _{i},
\\
\vartheta _{i,j} & := &
	\tau _{j-1,i}^- \tau _{j-1}^2 \tau _{j-1,i}
	=: \vartheta _{j,i}.
\end{eqnarray*}

\begin{statement}
For $k \geq 2$ let $B_k$ be generated by the set
$
\{ \tau _{i}^{(k)}; 1 \leq i \leq k-1 \}
$
according to Artin's relations.
The map
$
\tau _{i}^{(n)} \mapsto \tau _{i}^{(1+n)}
$
uniquely extends to a monomorphism
$
B_n \rightarrow B_{1+n}.
$
\end{statement}
The map extends to a homomorphism,
since the relations of $B_n$ are taken into relations
of $B_{1+n}$.
This homomorphism is injective,
since all additional relations of $B_{1+n}$ not being
relations of $B_n$ involve the generator $\tau _n$,
such that the
$\tau _i^{(1+n)}$ for $i \leq n-1$ obey the same relations
as do the generators of $B_n$.
\ep

\begin{statement}
Let
$
\{ t_{i}; 1\leq i \leq n-1 \}
$
be the set of permutations mapping
$ i \mapsto i+1, $
$ i+1 \mapsto i, $
$ \{ i, i+1 \} \not \ni j \mapsto j. $
The mapping
$
\pi :\tau _{i} 	\mapsto t_{i}
$
uniquely extends to an epimorphism
$
\pi \in \Hom (B_n, S_n)
$
from the braid group onto the group $S_n$ of
permutations of the numbers
$\{ 1, \ldots , n \}$
by mapping
$
\alpha \beta \mapsto \pi(\alpha ) \pi(\beta)
$.
\end{statement}
The $t_i$ exactly
satisfy
$t_{i}^{2} = 1$
and the Artin relations of the braid group.
Furthermore they generate $S_n$:
an $n$-permutation $\pi $ can be achieved by first
commuting the element
$\pi ^-(n)$ into its place by the $t_i$,
leading to a problem in $n-1$ elements.
\ep

\begin{statement}
Let $P_n$ be the kernel of
$ \pi \in \Hom (B_n, S_n) $,
called the pure or non-permuting braid group.
The elements
$$
\vartheta _{i,j} :=
	\tau _i^-\tau _{1+i}^-\ldots \tau_{j-2}^- \tau _{j-1}^2
	\tau _{j-2}\ldots \tau _{1+i}\tau _i
	=: \vartheta _{j,i},
$$
$i<j$, $i,j\in \{ 1,\ldots ,n\} $
generate $P_n$.
The relations between the generators
$ \vartheta _{i,j} $
are the consequences of
\begin{equation}
\label{PureBraids}
\vartheta _{k,l}^{-\epsilon } \vartheta _{i,j} \vartheta _{k,l}^{\epsilon}
=
\left\{ \begin{array}{ll}
	\vartheta _{i,j},
		& i < k \mbox{ or } l < i
\\
	\Ad ((\vartheta _{k,j} \vartheta _{l,j})^{\epsilon})
	(\vartheta _{i,j}),
		& i = k \mbox{ or } i = l
\\
	\Ad ([\vartheta _{l,j}^{\epsilon},\vartheta _{k,j}^{\epsilon}]
	^{-\epsilon })
        (\vartheta _{i,j}),
		& k < i < l
\end{array} \right. ,
\end{equation}
$i<j$, $k<l<j$, $\epsilon \in \{ -1, 1\} $,
$ \Ad (x)(y) := xyx^-,$
$ [x,y] := x y x^- y^-, $
for $x,y \in P_n $.
\end{statement}
A proof might proceed by finding a Schreier transversal
(cross-section)
to $B_n/P_n$ in $B_n$.
Cf.\cite{Birman}, pp.20,
where two pairs of parentheses are missing in lines $10$ and $11$.
A detailed algebraic proof is given in
\cite{Hansen}, pp.153.
\ep

The next theorem by recursion defines a
normal series for $P_n$.
\begin{statement}
\label{TheoremPureStructure}
Let
$
U_{n} \leq P_n \leq B_n
$
be the subgroup of the non-permuting braid group
generated by the set
$
\{ \vartheta _{n,i}; i \in \{ 1, \ldots , n-1 \} \}
$.
(This is the subgroup of the braid group $B_n$
in which only the $n$-th string is not held fixed.)
Then $P_n$ is isomorphic to the semidirect product
$
P_n \simeq P_{n-1} \opon U_{n},
$
where $P_{n-1}$ is imbedded into $P_n$
(by identifying the first
$n-1$ strings):
$ \vartheta _{i,j}^{(n-1)} \mapsto \vartheta _{i,j}^{(n)}. $
$P_{n-1}$ acts onto $U_{n}$ by conjugation.
\end{statement}
Statement and proof following
\cite{Birman}, pp.22.
By the relations eqn.~(\ref{PureBraids})
$U_n$ is normal in $P_n$.
The quotient $P_n / U_n$ is isomorphic to
the group generated by
$ \{ \vartheta _{i,j}; 1 \leq i < j \leq n-1 \} $,
that is $P_{n-1}$.
$P_{n-1}$ is imbedded into $P_n$ s.t. the exact sequence
$
1 \rightarrow
U_{n} \rightarrow
P_{n} \rightarrow
P_{n-1} \rightarrow
1
$
splits.
So the result follows,
cf.\cite{Robinson}, p.304.
\ep

Let
$
X_n := \{ (x_1,\ldots ,x_n) \in \C ^n;
\mbox{ if } i \neq j \mbox{ then } x_i \neq x_j \}
$
be the configuration space
of $n$ ordered distinct points in the plane.
For $x \in X_n$ let
$
Y_{x,m} :=
 \{ (y_1,\ldots ,y_m) \in X_m; y_i \neq x_j \}
$
be the configuration space of $m$ ordered distinct points
in the plane with $n$ points $x_j$ removed.

\begin{statement}
The projection
$
Y_{x,m+k} \rightarrow  Y_{x,m}
$
onto the first $m$ components
is a fiber bundle with
the fiber over
$y \in Y_{x,m}$
being
$
Y_{(x,y),k} =
 \{ (z_1,\ldots ,z_k) \in X_k; z_i \neq x_j, z_i \neq y_l \} .
$
\end{statement}
A theorem of Faddell and Neuwirth,
cf.\cite{Birman}, thm.1.2, pp.12.
Or cf.\cite{Hansen}.
\ep

\begin{statement}
The projection
$
X_n \rightarrow X_n/S_n
$
under the symmetric group $S_n$
is a regular $n!$ sheeted $S_n$ covering.
\end{statement}
Cf.\cite{Birman}, prop.1.1, p.11.
\ep

The following result on the fiber $Y_{x,m}$ is needed
for the algebraical approach to homology groups.
\begin{statement}
$Y_{x,m}$ is an Eilenberg-MacLane space of type
$(\pi _1(Y_{x,m}),1)$.
\end{statement}
For every $k$ we choose
$x \in X_k$
and set
$ Y_{k,m} = Y_{x,m}. $
The fiber bundle
$
Y_{n+m,k} \rightarrow Y_{n,m+k} \rightarrow Y_{n,m}
$
has the following exact homotopy sequence,
$$
\ldots
\stackrel{\partial}{\longrightarrow}
\pi _i(Y_{n+m,k})
\stackrel{i_*}{\longrightarrow}
\pi _i(Y_{n,m+k})
\stackrel{p_* \circ j}{\longrightarrow}
\pi _i(Y_{n,m})
\stackrel{\partial}{\longrightarrow}
\pi _{i-1}(Y_{n+m,k})
\ldots
$$
By induction assumption,
$
\pi _i(Y_{n+m,k}) \simeq \{ 1 \} \simeq \pi _i(Y_{n,m})
$
for every $i > 1$, every $n \geq 0$ and some $m$ and $k$.
This is true for $k = m = 1$ at least,
since in this case
$Y_{n+m,k}$ and $Y_{n,m}$ are Euclidean planes
with $n+m$ and $n$ punctures, respectively.
We conclude that also the homotopy groups
$ \pi _i(Y_{n,m+k}) $
vanish for $i > 1$.
\ep

We need further preparation in order to determine the structure
of the fundamental group $\pi _1(Y_{x,m})$.

\begin{statement}
\begin{enumerate}
\item
The normal subgroup
$ U_{n} \normal P_n $
is free over the set
$ \{ \vartheta _{n,i}; i \in \{ 1, \ldots , n-1 \} \} $.
\item
The non-permuting braid group $P_n$ is isomorphic
to the fundamental group $\pi _1(X_n)$.
\item
We have an isomorphism
$ B_n \simeq \pi _1(X_n / S_n). $
\end{enumerate}
\end{statement}
This topological proof follows
\cite{Birman}, pp.22.
For an algebraical proof of the first statement,
cf.\cite{Hansen}, pp.153.

1) From the homotopy sequence of the fiber bundle
$
Y_{n-1,1} \rightarrow X_n \rightarrow X_{n-1}
$
we obtain the exact sequence
\begin{equation}
\label{PureHomotopSequence}
1 \stackrel{\partial}{\longrightarrow}
\pi _1(Y_{n-1,1})
\stackrel{i_*}{\longrightarrow}
\pi _1(X_{n})
\stackrel{p_* \circ j}{\longrightarrow}
\pi _1(X_{n-1})
\stackrel{\partial}{\longrightarrow}
1.
\end{equation}
By topological reasons we can find
classes of paths
$ [ \gamma (\vartheta _{i,j}) ] \in \pi _1(X_k) $
corresponding to the generators of $P_k$
and obeying the
pure braid relations, at least.
So for
$ k \in \{ n-1 , n \} $
there are epimorphisms
$ h_k \in \Hom (P_k, \pi _1(X_k)) $
(which, in fact, will be shown to be isomorphisms,
once the present statement has been proved)
sending
$ \vartheta _{i,j} \mapsto [ \gamma (\vartheta _{i,j}) ] . $
The restriction
$ h_n \mid _{U_n} $
maps a generator
$ \vartheta _{n,j} $
to a loop running once around the $j$-th puncture,
starting and terminating at the basepoint
of the punctured plane
$Y_{n-1,1}$.
Since
$ \pi _1(Y_{n-1,1}) $
is a free group over these $n-1$ loops,
we have surjections
$
\pi _1(Y_{n-1,1}) \rightarrow
U_n \rightarrow
\pi _1(Y_{n-1,1}),
$
where the first map sends a free generator,
i.e.~a loop around the $j$-th puncture,
to $\vartheta _{n,j}$.
Since a free group of finite rank is Hopfian,
cf.\cite{Robinson}, thm.6.1.12, p.159,
this composed surjective endomorphism must be an automorphism
and $U_n$ must be isomorphic to
$
\pi _1(Y_{n-1,1}).
$

2) The exact sequence
$
1 \rightarrow U_n \rightarrow P_n \rightarrow P_{n-1}
\rightarrow 1
$
together with the
short exact sequence
(\ref{PureHomotopSequence})
of homotopy groups
yields a commutative diagram
where
$ U_n \rightarrow \pi _{1}(Y_{n-1,1}) $
is an isomorphism.
Assume,
$ P_{n-1} \rightarrow \pi _1(X_{n-1}) $
for some $n$ is an isomorphism as well.
Then also
$ P_n \rightarrow \pi _1(X_{n}) $
must be such.
Since
$ \pi _1(X_1) \simeq P_1 $
are trivial groups,
the statement is proved by induction.

3) This statement follows from the exact homotopy sequence
$
1 \rightarrow \pi _1(X_n) \rightarrow \pi _1(X_n / S_n)
\rightarrow S_n \rightarrow 1
$
of the regular covering
$ X_n \rightarrow X_n / S_n. $
One uses the fact that there are isomorphisms
$ P_n \rightarrow \pi _1(X_n) $
and a commutative diagram of two short exact sequences.
\ep

Therefore with the help of theorem
(\ref{TheoremPureStructure})
we may identify
$ P_n $
with the iterated semidirect product
$ F_2 \opon F_3 \opon \ldots F_{n-1} $
of free groups.

\begin{statement}
Let
$
t_{i,j} := t_i t_{1+i} \ldots t_{j-1} \in S_n
$
in the permutation group.
For every
$\pi \in S_n$
there is a sequence
$
\sigma _n(\pi ) = (t_{i_1,n}, t_{i_2,n-1}, \ldots , t_{i_n,1}),
$
the Schreier normal form of $ \pi $,
uniquely determined by the requirements that
$ 1 \leq i_k \leq 1 + n - k $
and that
$
\pi = t_{i_1,n} t_{i_2,n-1} \ldots t_{i_n,1}.
$
\end{statement}
An $n$-permutation $\pi $ can be achieved by first
commuting the element
$\pi ^-(n)$ into its place
(we are considering right-actions)
by applying
$t_{\pi ^-(n),n}$.
This leads to a problem in $n-1$ elements
\ep

\begin{statement}
For
$ j \in \{ 1, \ldots , n \} $
let
$ U_{j} \leq P_j \leq B_n $
be the subgroup of the non-permuting braid group $P_j$
generated by the set
$ \{ \vartheta _{j,i}; i \in \{ 1, \ldots , j - 1 \} \} $.
For every braid
$ \alpha \in B_n $
there is a sequence
$ (\tau _{i_1,n},\ldots , \tau _{i_{n-1},2}) $
with
$ i_j \in \{ 1, \ldots , n+1-j \} $
and a sequence
$ (\vartheta _{j_1}, \ldots , \vartheta _{j_{n-1}}) $
with
$ \vartheta _{j_k} \in U_{1+k} $
both uniquely determined by
$$
\alpha =
\tau _{i_1,n} \ldots \tau _{i_{n-1},2}
\vartheta _{j_1} \ldots \vartheta _{j_{n-1}} .
$$
The map sending a braid $\alpha $
to this pair of sequences is called Artin's
combed normal form of the braid group.
\end{statement}
We follow \cite{Birman}, cor.1.8.2, pp.24.
The set of products of the elements
$ \tau _{i_1,n}, \ldots , \tau _{i_{n-1},2} $
is a transversal (cross-section) of
$S_n$ in $B_n$ with respect to $P_n$.
This can be seen from the Schreier normal form for
the permutation group $S_n$ together with the homomorphism
$ \pi \in \Hom (B_n,S_n). $
It therefore exists such a unique product and a unique element
$ \vartheta \in P_n $
such that
$\alpha = \tau _{i_1,n} \ldots \tau _{i_{n-1},2} \vartheta. $
The structure of $P_n$ as an $n-1$-fold semidirect product
of free groups,
cf.thm.(\ref{TheoremPureStructure}),
yields a
unique decomposition of the element $\vartheta .$
\ep

\begin{statement}
\label{TheoremBraidSemiFree}
Let
$
P_{n}^{n} =
\{ \alpha \in B_n; \pi(\alpha)(n) = n \}
\leq B_n
$
be the subgroup of the braid group
that does not permute the $n$-th string
(the $n$-pure braid group).
Then $P_{n}^{n}$ is isomorphic to the semidirect product
$
P_{n}^{n} \simeq B_{n-1} \opon U_{n},
$
where $B_{n-1}$ is imbedded into $B_n$ by identifying the first
$n-1$ strings and it acts onto $U_{n}$ by conjugation.
\end{statement}
Let
$ \alpha = \tau _{i_1,n} \ldots \tau _{i_{n-1},2} \vartheta $
be given in the combed normal form.
If
$ \pi (\alpha)(n) = n, $
then $i_1 = n$ s.t.~$\tau _{i_1,n} = 1$.
So $ \alpha $ is an element of $B_{n-1}U_n$
and
$ P_n^n = B_{n-1}U_n. $
$U_n$ is normal in $P_n^n$ and
by the combed normal form,
$ B_{n-1} \cap U_n = \{ 1 \} ,$
so the claim follows.
For an algebraical proof of this decomposition,
cf.\cite{Chow1948} or \cite{Hansen}, pp.153.
\ep

In the following,
$F_n$ will be the free group generated by
$ \{ f_1, \ldots ,f_n \} . $
Subscripts may be dropped, if the rank of the group is
clear or immaterial.

The explicit formula
of the following imbedding theorem
allows the application of section
\ref{SectionMagnus}
to the braid group.

\begin{statement}[E.~Artin]
\label{TheoremArtinImbedding}
Let
$ \iota \in \Hom (U_{1+n},F_n) $,
$ \vartheta _{1+n,i} \mapsto f_i $
be an isomorphism of free groups.
Let
$
\psi \in \Hom (B_n, \Aut (F_n))
$
be such that
$
\psi (\alpha )(f) = \iota ( \alpha ^- \iota ^{-1}(f) \alpha )
$,
$ \alpha \in B_n, $
$ f \in F_n. $
Then $\psi $ is a monomorphism and
$$
\psi(\tau _{i}):f _{j} \mapsto
\left\{
       \begin{array}{cl}
       f _{i}f _{i+1}f _{i}^{-1},	&j=i,
\\
       f _{i},				&j=i+1,
\\
       f _{j},				&j\not \in \{ i,i+1\}
       \end{array}
\right. .
$$
\end{statement}
The existence of $\psi $ as a homomorphism
acting as claimed follows from the fact that
$U_{1+n} \simeq F_n$ is invariant under conjugation
with elements from $B_n$
and from the
relations between the elements
$\tau _i$ and $\vartheta _{1+n,j}$.
For a proof that the kernel of $\psi $
is trivial,
cf.\cite{Birman}, cor.1.8.3, pp.25.
\ep

\begin{statement}
\label{TheoremFreeGroupAction}
The action of the elements
$
\vartheta _{n,i} \in B_{n}
$
generating the free group
$U_n$ of rank $n-1$ in $B_n$
onto the generators of the free group
$ F_{n} $
is given by
($ \epsilon \in \{-1,1\} $)
$$
\vartheta _{n,i}^{\epsilon }(f_j) =
\left\{ \begin{array}{ll}
f_j,
	&j<i
\\
\Ad ((f_if_n)^{\epsilon })(f_j),
	&j \in \{ i, n \}
\\
\Ad ([f_{n}^{\epsilon },f_i^{\epsilon }]^{-\epsilon })(f_j),
	&i<j<n
\end{array} \right. ,
$$
where
$ \Ad (x)(y) := xyx^-, $
$ [x, y] := xyx^-y^- $ for $x,y \in F_n$.
\end{statement}
This can be inferred from equation
(\ref{PureBraids})
by imbedding
$ U_n \leq P_n \rightarrow P_{1+n} $
and
$ F_n \simeq U_{1+n} \leq P_{1+n} . $
\ep

The procedure of taking a semidirect product
$ B_n \opon F_n $
and,
by theorem (\ref{TheoremBraidSemiFree}),
imbedding it into the braid group $B_{1+n}$,
$ B_n \opon F_n \simeq P_{1+n}^{1+n} \leq B_{1+n}, $
can be iterated.
We set
\begin{eqnarray*}
B_{n,m}	& = & B_n \opon F_n \opon \ldots \opon F_{n+m-1},
\\
F_{n,m}	& = & F_n \opon \ldots \opon F_{n+m-1},
\end{eqnarray*}
such that
$ B_{n,m} = B_n \opon F_{n,m} \leq B_{n+m}. $
In this iterated product,
$F_k$ is acting onto $F_{1+k}$ according to
thm.(\ref{TheoremFreeGroupAction}).
$F_{n,m}$ may be identified
with the subgroup of the pure braid group
$P_{n+m},$
in which only the last $m$ strings are
not held fixed.

\begin{statement}
The fundamental group of the space
$
Y_{x,m} =
 \{ (y_1,\ldots ,y_m) \in X_m; y_i \neq x_j \}
$
for $x \in X_n$
is isomorphic to the group $F_{n,m}$.
\end{statement}
The fiber bundle
$
Y_{n,m} \rightarrow X_{n+m} \rightarrow X_{n}
$
due to vanishing of the higher homotopy groups
yields the exact sequence
$
1 \rightarrow
\pi _1(Y_{n,m}) \rightarrow
\pi _1(X_{n+m}) \rightarrow
\pi _1(X_{n}) \rightarrow
1.
$
We know,
$
\pi _1(X_{n+m}) \simeq P_{n+m}
$
and
$
\pi _1(X_{n}) \simeq P_n.
$
The exact sequence
$
1 \rightarrow
F_{n,m} \rightarrow
P_{n+m} \rightarrow
P_n \rightarrow
1
$
can be derived from the combed normal form for
$ P_{n+m}. $
By the five lemma,
$ F_{n,m} \simeq \pi _1(Y_{n,m}). $
\ep

\begin{statement}
The fiber
$
Y_{x,m} =
\{ (y_1, \ldots ,y_m) \in X_m; y_i \neq x_j \}
$
over
$ x \in X_n $
has homology
$
H_*(Y_{x,m},M) \simeq H_*(F_{n,m},M)
$
for any $F_{n,m}$ module $M$.
\end{statement}
We have shown,
$Y_{x,m}$ is an Eilenberg-MacLane space of type
$(F_{n,m},1)$.
As a smooth manifold,
it is homotopic to a cell complex,
cf.\cite{Milnor}, thm.3.5, pp.20 and p.36.
So the claim follows e.g.~from
\cite{BrownCohomology}, par.III.1,
pp.56, using Eilenberg's theorem in
\cite{Whitehead}.
\ep

Now we are in the position to show,
in which way the representations of $B_n$ on
$
H_*(Y_{x,m},\chi )
$
derived in
\cite{Lawrence1990}
can be constructed in purely algebraic terms.

\section{%
	Resolution of the integers over $ F_{n,m} $
}
We construct a free
$ F_{n,m} \simeq \pi _1(Y_{x,m}) $
resolution $C$ of the integers.
It is chosen in such a way
that there is a particular braid action onto complexes
$ N \otimes _{F_{n,m}} C $
(and onto their homology)
for any right $B_{n,m}$ module $N$.
This action is described by
the braid valued Burau matrices
that will be derived in the next section.

For
fixed integers $n, m$
and for
$ k \in \{ 1, \ldots , m \} $
let
$ F_{n+k-1} = F^k $
be the free group freely generated by the subset
$
\{ f_i^{(k)}; i \in \{ 1, \ldots , n+k-1 \} \} \subset F^k
$.
Let $I^k$ be the relative augmentation ideal of
$F^{k}$
in the group ring of
$
F_{n,k} = F^{1} \opon \ldots \opon F^{k}
$.
$I^k$
is freely generated as a left module over
$
\{ (f_i^{(k)} - 1); i \in \{ 1, \ldots , n+k-1 \} \}
$.
We distinguish $I^k$
from its image in $\Z F_{n,k}$.
Therefore we write
$
s_i^k
$
when we consider the generating elements in the ideal
$I^k$ rather than in $\Z F_{n,k}$.
The natural imbedding
$\iota : I^k \rightarrow \Z F_{n,k} $
thus sends
$s_i^k \mapsto (f_i^{(k)} - 1). $
We consider the ideal $ I^k $ as an $F_{n,k}$  bimodule
(and therefore as a bimodule over
$
F_{n,j} \leq F_{n,k},
$
$
j \leq k
$).
We sometimes write $s^k$ for
$
(f^{(k)} - 1) \in I^k,
$
with
$ f^{(k)} \in F^k. $

\begin{statement}
\label{TheoremF_{n,m}Complex}
Let
$$
C_k^m :=
\oplus _i \Z F_{n,m} \otimes _{F_{n,i_1}} I^{i_1}
	\otimes _{F_{n,i_2}} \ldots \otimes _{F_{n,i_k}} I^{i_k}
$$
for
$
k \in \{ 1, \ldots , m \} ,
$
$
1 \leq i_k < i_{k-1} < \ldots < i_1 \leq m
$
and let
$
C_0^m := \Z F_{n,m}
$
with the natural structures as left
$F_{n,m}$ modules.
Let
$
\partial _k^m : C_k^m \rightarrow C_{k-1}^m
$
be given by
$
\partial _k^m  :=  \sum _{l=1}^k (-)^{1+l} \partial _{k,l}^m,
$
with
\begin{eqnarray*}
\lefteqn{
\partial _{k,l}^m :
	 \Z F_{n,m} \otimes _{F_{n,i_1}} I^{i_1}
	\ldots
	\otimes _{F_{n,i_{l-1}}}
	I^{i_{l-1}} \otimes _{F_{n,i_{l}}} I^{i_l}
	\otimes _{F_{n,i_{l+1}}}
	\ldots I^{i_k}
}
& & \\
& &
\rightarrow
	\Z F_{n,m}
	\otimes _{F_{n,i_1}} I^{i_1}
	\ldots
	\otimes _{F_{n,i_{l-1}}} I^{i_{l-1}} \iota(I^{i_l})
	\otimes _{F_{n,i_{l+1}}}
	\ldots I^{i_k},
\end{eqnarray*}
where
$
\partial _{k,l}^m :=
1 \otimes \ldots \otimes \iota \otimes \ldots \otimes 1
$
and the imbedding
$ \iota : I^l \rightarrow \Z F_{n,l} $
occurs at the $(l+1)$-st position in the $(k+1)$-fold
tensor product.
Then
$$
0
\rightarrow
C_m^m
\stackrel{\partial _m^m}{\longrightarrow}
C_{m-1}^m
\stackrel{\partial _{m-1}^m}{\longrightarrow}
\ldots
\stackrel{\partial _1^m}{\longrightarrow}
C_{0}^m
\stackrel{\epsilon}{\longrightarrow}
\Z
\rightarrow
0
$$
with the augmentation
$
\epsilon : \Z F_{n,m} \rightarrow \Z
$
sending
$F_{n,m}$ to $1$
is a complex of free left $F_{n,m}$ modules,
augmented over the trivial left $F_{n,m}$ module $\Z $.
\end{statement}
We will fix $m$ and drop corresponding superscripts.
Each $C_k$ is a free left $F_{n,m}$ module,
since all the factors in the defining tensor products
are free over the respective rings.
$\partial _k$ is a homomorphism of
left $F_{n,m}$ modules and
$\Image (\partial _1)$
is contained in the augmentation
ideal of $F_{n,m}$ in $\Z F_{n,m},$ s.t.
$
\Image (\partial _1) \leq \Ker (\epsilon ).
$
Furthermore for
$
k \in \{ 2, \ldots ,m \} ,
$
$
\Image (\partial _k) \leq \Ker (\partial _{k-1}),
$
i.e.~$ \partial _{k-1} \circ \partial _k = 0, $
since in the sum expressing
$
\partial _{k-1} \circ \partial _k,
$
each tensor product with
$\iota $ at the $i$-th and at the $j$-th position,
$i \neq j ,$
occurs twice,
once with sign
$(-)^i(-)^{j-1}$
and once with sign
$(-)^j(-)^i.$
\ep

In explicit terms, the components of the
boundary operator act as
$$
\partial _{k,l}^m(
s^{i_1} \otimes \ldots
s^{i_{l-1}} \otimes s^{i_l}
\otimes \ldots  s^{i_k} )
 =  s^{i_1} \otimes \ldots
s^{i_{l-1}}(f^{(i_l)} - 1) \otimes
\ldots s^{i_k},
$$
where the product
$
s^{i_{l-1}}(f^{(i_l)} - 1)
$
has to be understood
as the right action of
$
I^{i_l} \leq \Z F_{n,i_l}
$
onto the
bimodule $I^{i_{l-1}}$.
This expression for the boundary
might be compared with
the topological boundary matrix in
\cite{Lawrence1990}, sec.3.2, p.152.

The complexes
$
(C^m,\partial ^m)
$
can be described in recursive terms:
for $k > 0$, $m > 0$,
the modules are given by
$
C_0^m = \Z F_{n,m},
$
$
C_1^1 \simeq I^1,
$
$
C_{1+k}^1 = \{ 0 \},
$
$$
C_k^{1+m} \simeq I^{1+m} \otimes _{F_{n,m}} C_{k-1}^m \oplus
	\Z F_{n,1+m} \otimes _{F_{n,m}} C_k^m.
$$
The homomorphisms are determined by
the augmentation
$
\epsilon ^m : C_0^m \rightarrow \Z ,
$
by
$
\partial _1^m : C_1^m \rightarrow C_0^m,
$
which imbeds the direct summands
$
I^k \rightarrow \Z F_{n,k}
$
for all
$ k \in \{ 1, \ldots , m \} $
and by the recursive step
\begin{eqnarray*}
\lefteqn{
\partial _k^{1+m} (j \otimes c + f \otimes d)
}
&&
\\
& = &
(\partial _1^{1+m}(j) \otimes c
- j \otimes \partial _{k-1}^m(c))
+ f \otimes \partial _k^m(d),
\end{eqnarray*}
for
$0 < k \leq m,$
$j \in I^{1+m},$
$c \in C_{k-1}^m,$
$f \in F_{n,1+m},$
$d \in C_{k}^m$.

\begin{statement}
The complex
$
0
\rightarrow
C_m^m
\stackrel{\partial _m^m}{\longrightarrow}
C_{m-1}^m
\stackrel{\partial _{m-1}^m}{\longrightarrow}
\ldots
\stackrel{\partial _1^m}{\longrightarrow}
C_{0}^m
\stackrel{\epsilon}{\longrightarrow}
\Z
\rightarrow
0
$
of left $F_{n,m}$ modules
is exact.
\end{statement}
Inductively for every $m > 0$
we will construct a free left $F_{n,m}$ resolution
of the integers that coincides with the
augmented complex $C^m$.
Clearly, for $m = 1$ we have the resolution
$
0 \rightarrow I_{F_{n,1}} \rightarrow \Z F_{n,1}
	\rightarrow \Z \rightarrow 0.
$
Now we assume, there is an $m$ s.t.
$
C^m \rightarrow \Z \rightarrow 0
$
is a free resolution and the short sequences
$
0 \rightarrow K_i^m \rightarrow C_i^m  \rightarrow K_{i-1}^m
	\rightarrow 0,
$
with
$
K_i^m = \Ker (\partial _i^m)
$
are exact.
In order to obtain the resolution $C^{1+m}$,
in the first step we construct
the diagram
$$
\begin{CD}
I^{1+m} @>>> K_0^{1+m} @>>> K_0^m
\\
@V{\Id}VV @VVV @VVV
\\
I^{1+m} @>>> C_0^{1+m} @>{p_0}>> C_0^m
\\
@V0VV @V{\epsilon ^{1+m}}VV @VV{\epsilon ^{m}}V
\\
0 @>>> \Z @= \Z
\end{CD}.
$$
This will be shown to be commutative
with exact rows and columns
where the second maps in each row and and in each column
are onto and the first maps are injective.
The map
$ \Z \rightarrow \Z $
is the identity.
The vertical map on the left,
$
I^{1+m} \rightarrow 0,
$
is obvious,
the vertical map on the right,
$
\epsilon ^m : C_0^m \rightarrow \Z ,
$
is the augmentation.
We set
$
C_0^{1+m} = \Z F_{n,1+m},
$
which is a free left $F_{n,1+m}$ module
and as an abelian group is isomorphic to
the direct sum
$
I^{1+m} \oplus C_0^m
$.
The map
$
p_0 : C_0^{1+m} \rightarrow C_0^m
$
is the surjection induced by
the quotient map modulo the normal subgroup
$
F^{1+m} \normal F_{n,1+m}
$
and so it has the kernel
$ I^{1+m}. $
The middle vertical map
$
\epsilon ^{1+m} : C_0^{1+m} \rightarrow \Z
$
is defined by commutativity of the
two lower squares.
It is the augmentation
$
\Z F_{n,1+m} \rightarrow \Z
$.
Therefore,
the kernels of the vertical maps are
$I^{1+m}$,
$
K_0^{1+m} = I_{F_{n,1+m}}
$
and
$
K_0^{m} = I_{F_{n,m}},
$
respectively.
The maps in the first row are the induced ones
and by the snake lemma,
cf.\cite{HiltonStammbach}, lemma III.5.1, p.99,
this row is exact
with the first map injective and the second one onto.
Now, since $I^{1+m}$ and
$
D_1 = \Z F_{n,1+m} \otimes _{F_{n,m}} C_1^m
$
are free left $F_{n,1+m}$ modules,
we obtain free presentations
$ 0 \rightarrow  I^{1+m} \rightarrow I^{1+m} $
and
$ L_1 \rightarrow  D_1 \rightarrow K_0^m $
of the kernels,
fitting into the following diagram
$$
\begin{CD}
0 @>>> K_1^{1+m} @>{p_1}>> L_1
\\
@VVV @VVV @VVV
\\
I^{1+m} @>>> C_1^{1+m} @>{p_1}>> D_1
\\
@V{\Id }VV @V{\partial _1^{1+m}}VV @VV{\Delta _1}V
\\
I^{1+m} @>>> K_0^{1+m}  @>{p_0}>> K_0^{m}
\end{CD}.
$$
Again we will choose maps in such a way
that exactness of all sequences is preserved,
if we attach zero maps
on both sides of the sequences
and that the diagram is commutative.
The vertical map on the left,
$
I^{1+m} \rightarrow I^{1+m},
$
is the identity,
the vertical map on the right,
$
\Delta _1 : D_1 \rightarrow K_0^m,
$
is
$
\Delta _1 = \eta \otimes \partial _1^m
$,
$ \eta : F_{n,1+m} \rightarrow F_{n,m} $
being the augmentation homomorphism
with respect to $F^{1+m}$.
$
C_1^{1+m} = I^{1+m} \oplus D_1
$
as free left $F_{n,1+m}$ modules
and the middle row is
induced by this decomposition.
The free presentation
$
\partial _1^{1+m} : C_1^{1+m} \rightarrow K_0^{1+m}
$
we are on the way of constructing
as the middle column must make the lower
squares commute.
It is determined by its components.
$
\partial _{1,1}^{1+m}: I^{1+m} \rightarrow K_0^{1+m}
$
is the imbedding
$ s_i^{1+m} \mapsto (f_i^{(1+m)} - 1) $
and
$
\partial _{1,2}^{1+m}:
D_1 \rightarrow K_0^{1+m}
$
maps
$
f \otimes _{F_{n,m}} s_{i}^{k} \mapsto f (f_i^{(k)} - 1)
\in I_{F_{n,1+m}} = K_0^{1+m}
$.
Commutativity of
the lower left square is obvious.
To show
$
\Delta _1 = p_0 \circ \partial _{1,2}^{1+m},
$
we notice,
$
p_0(\partial _{1,2}^m (f g \otimes s_i^k))
=
f (f_i^{(k)} - 1)
$
for
$ f \in F_{n,m}, $
$ g \in F^{1+m}. $
On the other hand,
$
\Delta _1(fg \otimes s_i^k)
=
f (f_i^{(k)} - 1).
$
$\Delta _1$ is onto $K_0^m$,
since
$\partial _1^m$ is,
by assumption.
So by the five lemma also
$
\partial _1^{1+m}
$
is onto.
The kernel of the left vertical map
is $0$.
The kernel of $\Delta _1$
can be seen to be the module
$
L_1 = (I^{1+m} \otimes _{F_{n,m}} C_1^m) \oplus K_1^m
$,
using the decomposition
$
\Z F_{n,1+m} \simeq I^{1+m} \oplus \Z F_{n,m},
$
where both direct sums are to be understood
as direct sums of abelian groups,
rather than $F_{n,1+m}$ modules.
Again by the snake lemma, the row of kernels
is exact, with a surjection after an injection.
Once more we will repeat the construction of
free presentations of the kernels occuring
in the top row.
This leads to an iterative procedure which
eventually reaches an isomorphism s.t.~the
resolutions are finite,
i.e.~the vertical sequences begin with $0$ at the top.
Thus we consider the diagram
$$
\begin{CD}
0 @>>> K_2^{1+m} @>>> L_2
\\
@VVV @VVV @VVV
\\
0 @>>> C_2^{1+m} @= D_2
\\
@VVV  @V{\partial _2^{1+m}}VV @VV{\Delta _2}V
\\
0 @>>> K_1^{1+m} @>{p_1}>> L_1
\end{CD}.
$$
Here
$
C_2^{1+m} = D_2 =
I^{1+m} \otimes _{F_{n,m}} C_1^m \oplus
\Z F_{n,1+m} \otimes _{F_{n,m}} C_2^m.
$
The map
$
D_2 \rightarrow L_1
$
is
$
\Delta _2 =
\partial _1^{1+m} \otimes 1 + 1 \otimes \partial _2^m.
$
It is onto,
since by assumption
$\partial _2^m$
is onto $K_1^m$,
and
$\partial _1^{1+m}$
is onto $I^{1+m}$.
The free presentation
$
\partial _2^{1+m} : C_2^{1+m} \rightarrow K_1^{1+m}
$
as a $F_{n,1+m}$ module map is uniquely determined
by the requirement for commutativity
of the lower right square,
$
p_1 \circ \partial _2^{1+m} = \Delta _2.
$
The components of
$
\partial _2^{1+m}
$
are
$
( \partial _1^{1+m} \otimes 1 - 1 \otimes \partial _1^m )
$
and
$
1 \otimes \partial _2^m
$
and it explicitly acts as
\begin{eqnarray*}
\lefteqn{%
	\partial _2^{1+m}(
	s_i^{1+m} \otimes _{F_{n,m}} c_1
	\oplus
	f^{(1+m)} \otimes _{F_{n,m}} c_2
	)
}
&&
\\
& = &
(f_i^{(1+m)} - 1) \otimes c_1
- s_i^{1+m} \partial _1^m(c_1)
+ f_i^{(1+m)} \otimes _{F_{n,m}} \partial _2^m(c_2),
\end{eqnarray*}
for
$c_1 \in C_1^m$, $c_2 \in C_2^m$,
$f^{(1+m)} \in F_{n,1+m}$.
The map
$\partial _2^{1+m}$
is onto,
since $\Delta _2$ is.
The kernel of the map
$ \Delta _2 $
is the module
$
L_2 =
\Z F_{n,1+m} \otimes K_2^m +
< s_i^{1+m} \otimes \partial _2^m(c)
- (f_i^{(1+m)} - 1) \otimes c>.
$
To validate this, one uses the decomposition
$ \Z F_{n,1+m} \simeq I^{1+m} \oplus \Z F_{n,m} $
and the linear independence of
$ \{ f_i^{(1+m)} - 1 \} . $
The snake lemma again shows that the top
row is exact with a surjection after an
injective map.
Therefore
$
\Ker (\partial _2^{1+m}) = K_2^{1+m} = L_2.
$
Proceeding iteratively in this way,
using the induction assumption
we find,
$
C_{1+l}^{1+m} =
I^{1+m} \otimes _{F_{n,m}} C_l^m \oplus
\Z F_{n,1+m} \otimes _{F_{n,m}} C_{1+l}^m,
$
$
K_{1+l}^{1+m} =
\Z F_{n,1+m} \otimes K_l^m +
< s_i^{1+m} \otimes \partial _l^m(c)
- (f_i^{(1+m)} - 1) \otimes c>.
$
\ep

For every right $B_{n,m}$ module $N$,
the tensor complex
$
N \otimes _{F_{n,m}} C
$
and its homology
$ H_*(F_{n,m}, N) $
have a structure as
a right $B_n$ module.
The right action is given by
$
\alpha  :
n \otimes _{F_{n,m}} c
\mapsto
n \alpha \otimes _{F_{n,m}} \alpha (c),
$
where $\alpha \in B_n$ acts via Artin's automorphisms onto
$c \in C^k$.
We prefer a different point of view.
Instead of the augmentation ideals
$ I^l = I_{F_{n+l-1},F_{n,l}} $
used in the construction of the resolution $C$,
we may consider the ideals
$ \bar{I}^l = I_{F_{n+l-1},B_{n,l}}, $
with
$ B_{n,l} = B_n \opon F_{n,l}. $
$\bar{I}^l$ as a left $B_{n,l}$ module is isomorphic to
$ \Z B_{n,l} \otimes _{F_{n,l}} I^l $
and as an ideal in the ring
$\Z B_{n,l}$,
it carries a structure as right $B_{n,l}$ module, too.
Let the complex $\bar{C}$ be constructed similar
to $C$ in thm.~(\ref{TheoremF_{n,m}Complex}) by setting
$$
\bar{C}_k^m :=
\oplus _i \Z B_{n,m} \otimes _{B_{n,i_1}} \bar{I}^{i_1}
	\otimes _{B_{n,i_2}} \ldots \otimes _{B_{n,i_k}} \bar{I}^{i_k}.
$$
This complex is not a free resolution of $\Z $ anymore
but it is still exact and
for any right $B_{n,m}$ module $N$ we may identify
$ N \otimes _{F_{n,m}} C \simeq N \otimes _{B_{n,m}} \bar{C} $
as right $B_n$ modules,
where the $B_n$ action on
$ N \otimes _{B_{n,m}} \bar{C} $
is more obvious.
This action is known
as soon as the representations
on the ideals
$ I_{F_k,B_kF_k} $
for all $k$
are known,
because
$ \bar{I}^l = I_{F_{n+l-1},B_{n,l}} $
is imbedded into
$ I_{F_{k},B_{k}F_{k}} $
with
$ k = n+l-1. $
This action
will be found in the next section.

\section{%
	The braid valued Burau matrices
        }
\label{SectionBurau}
The relative augmentation ideal
$I_{F_n,B_nF_n}$
is a right $B_n$ module
and is free
as a left $B_nF_n$ module.
The right $B_n$ action therefore is determined
by matrices
which also determine the action onto the complex
$ N \otimes _{F_{n,m}} C $
and onto its homology,
for any right $B_{n,m}$ module $N$.
They generalise the
Burau matrices.

The Burau representation of the braid group
is a ``classical''
Birman-Magnus representation in the sense
of section \ref{SectionMagnus},
cf.\cite{Birman,Magnus1974}.
Let
$ \psi \in \Hom (F_n,F_1=<t>). $
This map has the kernel
$
U =
\{ f_{i_1}^{\epsilon _1}\ldots  f_{i_k}^{\epsilon _k};
   \sum \epsilon _l = 0 \} .
$
According to thm.~(\ref{TheoremArtinImbedding}),
for every braid $ \alpha \in B_n $ and every $ f \in F_n $,
$ \psi (\alpha (f)) = \psi (f). $
So there is a Birman module
$ M_{\psi } \simeq \Z [t] \otimes _{B_nF_n} I_{F_n,B_nF_n} $
carrying a right action of $B_n$,
\begin{eqnarray*}
\begin{array}{rcl}
s_j \tau _{i} & = &
	\delta _{\psi }(\tau _{i}(f_{j}))
\\
	& = &
	\left\{ \begin{array}{rcll}
		\delta _{\psi }(f _{i}f _{i+1}f _{i}^{-1}) & = &
			(1 - t) s_{i} + t s_{i+1}, & j=i
\\
		\delta _{\psi }(f _{i}) & = & s_{i}, & j = i + 1
\\
		\delta _{\psi }(f _{j}) & = & s_{j}, & j \not \in \{ i,i+1 \}
	\end{array}\right. .
\end{array}
\end{eqnarray*}
This means the right action is determined by the well-known
Burau matrices
$$
\tau _i \mapsto
\left(
\begin{array}{ccccc}
{\mathbf 1}_{i-1}&0		&0	&0
\\
0	&( 1 - t )	&t	&0
\\
0	&1		&0	&0
\\
0	&0		&0	&{\mathbf 1}_{n-i-1}
\end{array}
\right) .
$$

The Magnus submodule
$
\delta _{\psi }(U) \simeq U / [U,U] \simeq H_1(F_n, \Z [t] )
$
according to thm.~(\ref{TheoremFoxImage}) is
$$
\{ \sum  g_j s_j;  g_j \in \Z [t], \sum g_j (\psi(f_j) - 1) = 0 \} =
 \{ \sum  g_j s_j; g_j \in \Z [t], \sum g_j = 0 \} ,
$$
since $ \psi (f_j) = t $
and the polynomial ring
$ \Z [t] $ is free of zero divisors.
%

The modules constructed in
section \ref{SectionMagnus}
according to
thms.~(\ref{TheoremSequence}) and (\ref{TheoremIdentify})
are induced by the representation on
the relative augmentation ideal
$ I_{F,EF}. $
In the last section we have seen
that the right $B_n$ action
onto
$ N \otimes _{F_{n,m}} C $
for any right $B_{n,m}$ module
also is determined
by the ideals
$ I_{F,BF}. $
So we obtain the
braid valued Burau representation,
\cite{ConstantinescuLuedde1992,Luedde1992}.
\begin{statement}
Let $M$ be a free left $BF$ module of rank $n$
with basis
$ \{ s_1, \ldots , s_n \} . $
Then the map
$
\rho : \{ \tau _1, \ldots , \tau _{n-1} \} \rightarrow
\End _{BF}(M)
$
defined by
$$
\rho (\tau _{i})(s_j) =
\tau _{i}
\left\{ \begin{array}{ll}
(1 - f_{i}f_{i+1}f_{i}^{-1})s_{i} + f_{i}s_{i+1}, & j=i
\\
s_{i},						& j = i + 1
\\
s_{j},						& j \not \in \{ i, i+1 \}
\end{array} \right.
,$$
uniquely extends to a monomorphism
$ \rho \in \Hom (B,\End _{BF}M) $.
\end{statement}
$I_{F,BF}$ is a free left $BF$ module over the set
$ \{ (f_i - 1); i \in \{ 1, \ldots , n \} \} ,$
so we identify $M$ with $I_{F,BF}$ via
$ s_i \mapsto (f_i - 1). $
By multiplication from the right
$I_{F,BF}$ is a module over $B$,
$
(f_j - 1) \mapsto
(f_j - 1) \tau _i = \tau _i (\tau _i(f_j) - 1).
$
The element $\tau _i(f_j)$ is determined by Artin's
thm.~(\ref{TheoremArtinImbedding})
so we obtain the equation as claimed.
\ep

In terms of matrices over the ring $\Z (BF)$
the representation of $\tau _i$ is given by
\begin{equation}
\label{BraidValuedBurau}
\tau _i \mapsto
\tau _i
\left(
\begin{array}{ccccc}
{\mathbf 1}_{i-1}&0				&0	&0
\\
0	&(1-f_{i}f_{i+1}f_{i}^{-1})	&f_{i}	&0
\\
0	&1		&0			&0
\\
0	&0		&0		&{\mathbf 1}_{n-i-1}
\end{array}
\right) .
\end{equation}

With the help of the imbedding
$
B_nF_n \rightarrow B_{1+n},
$
$
\tau _i f_j \mapsto \tau_i \vartheta _{1+n, j},
$
the semidirect product $B_nF_n$
and therefore the matrix elements of
eqn.~(\ref{BraidValuedBurau}) can again
be represented by the
braid valued Burau matrices.
By $m$-fold iteration
we obtain representations of $B_n$
(or of subgroups $ S \leq B_n $ )
in terms of matrices with values in the ring
$\Z (B_n F_n F_{1+n}\ldots F_{n+m-1})$
($ \Z (S F_n F_{1+n} \ldots F_{n+m-1}), $ respectively).

Finally we show
that the local coefficient system
chosen in \cite{Lawrence1990}
enables us to identify the recursion matrices
given there
with images of the braid valued Burau matrices.
We imbed
$ F_{n,m} = F_n \opon \ldots F_{n+m-1} $
into $B_{n+m}$.
The generators
$ \{ f_i^{(k)} ; i \in \{ 1, \ldots , n+k-1 \} \} $
of $F_{n+k-1}$ are mapped as
$ f_i^{(k)} \mapsto \vartheta _{k+n,i}. $
\begin{statement}
\label{TheoremLawrenceCoefficients}
Let a homomorphism
$
\chi \in Hom(F_{n,m}, \C \backslash \{ 0 \} )
$
be given by the map
$
\vartheta _{l,i} \mapsto q_{l,i},
$
for
$ l \in \{ 1+n, \ldots , m+n-1 \} , $
$ i < l . $
If
$
q_{l, i} = q_{l, j}
$
for
$ i, j \in \{ 1, \ldots , n \}, $
then $ \chi $ extends to a map
$
\bar{\chi } \in Hom(B_{n} \opon F_{n,m}, \C \backslash \{ 0 \} )
$
by setting
$
\bar{\chi } (\alpha \vartheta) = \chi (\vartheta),
$
$ \alpha \in B_n, $
$ \vartheta \in F_{n,m}. $
\end{statement}
Eqn.~(\ref{PureBraids})
shows that any map
$
\vartheta _{l,i} \mapsto q_{l,i}
$
of the generators into an abelian group
extends to a homomorphism of $F_{n,m}$.
This is also the content of lemma 2.1, p.145, \cite{Lawrence1990}.
The relations defining the semidirect product
$ B_{n,m} = B_n \opon F_{n,m} $
are those in $B_n$, in $F_{n,m}$ and in addition
those between $B_n$ and $F_{n,m}$:
$$
\tau _{i}^{-1} \vartheta _{l,j} \tau _i =
\left\{
       \begin{array}{cl}
       \vartheta _{l,i} \vartheta _{l,i+1} \vartheta _{l,i}^{-1},
	& j=i
\\
       \vartheta _{l,i},
	& j=i+1
\\
       \vartheta _{l,j},
	& j \not \in \{ i,i+1\}
       \end{array}
\right. .
$$
So the conditions $q_{l,i} = q_{l,j}$ allow
the extension of $\chi $ to $B_{n,m}$ by mapping
$B_n$ to $1$.
\ep

{}From the following matrices,
representations of $B_n$
on the homology $H_m(Y_{x,m},\chi )$
are derived
in \cite{Lawrence1990}.
\begin{statement}[R.J.Lawrence]
For
$ r \in \{ 0, 1, \ldots , m \} , $
$ l \in \{ 1 + r, \ldots , m + n - 1 \} , $
let the matrices $ A_l^{(r)} $
be recursively defined by the equations
\begin{eqnarray*}
A_i^{(0)} & = & 1,
\\
A_i^{(k)} & = &
\left(
\begin{array}{cccc}
{\mathbf 1}_{m + n - k - i} & 0 & 0 & 0
\\
0 & 0 & b_{k,1+i}^{(k-1)} & 0
\\
0 & 1 & (1 - (b_{k,1+i}^{(k-1)})^{-1}
		b_{k,i}^{(k-1)}
		b_{k,1+i}^{(k-1)})
			& 0
\\
0 & 0 & 0 & {\mathbf 1}_{i - 2}
\end{array}
\right) A_i^{(k-1)},
\end{eqnarray*}
where for
$ q_{1+r,p} \in \C \backslash \{ 0 \} , $
$ p > 1+r $
we set
\begin{eqnarray*}
b_{1+r,p}^{(r)} & = & A_{1+r,p}^{(r)2} / q_{1+r,p},
\\
A_{1+r,p}^{(r)} & = &
	A_{1+r}^{(r)} A_{2+r}^{(r)} \ldots A_{p-1}^{(r)}
	A_{p-2}^{(r)-} \ldots A_{1+r}^{(r)-}.
\end{eqnarray*}
If
$
q_{j, s} = q_{j, t}
$
for
$ j \in \{ 1, \ldots , m \} , $
$ s, t \in \{ 1 + m, \ldots , n + m \} , $
then the matrices
$ A_{l}^{(m)} , $
$ l \in \{ 1 + m, \ldots , m + n - 1 \} $
generate a homomorphic image of the braid group $B_n$.
\end{statement}
This is thm.~3.4, pp.156, \cite{Lawrence1990},
where it is proved
by computing the braid action on a cell complex.
We have only given the recursion matrices
$ A_l^{(r)} = A_{l,1+l}^{(r)} $
instead of the more complicated ones
$ A_{l,p}^{(r)}, $
because the matrices
$ A_l^{(r)} $
generate the others.
We have furthermore used lemma 3.3, p.153,
parts (i) and (ii) to slightly change the form
of $ A_l^{(r)} $
compared to the expression in loc.~cit.
\ep

The theorem can be understood in our algebraic approach.
We compare
$ A_l^{(r)} $
to the braid valued Burau matrix,
eqn.~(\ref{BraidValuedBurau}).
Lawrence's matrices act by left multiplication with column vectors,
the generalised Burau matrices act
by right multiplication with row vectors.
Furthermore, the conventions
in \cite{Lawrence1990}
on the imbedding
$B_{n,m} \rightarrow B_{n+m}$
are different from ours.
This accounts for the reversal of products of matrices,
for rearranging of rows and columns
and for reindexing of the generators.
Apart from this both matrices have the same appearance.
The lowest order matrices due to lemma
(\ref{TheoremLawrenceCoefficients})
yield a representation
$ \bar{\chi } \in \Hom (B_{n,m}, \C \backslash \{ 0 \} ), $
$ \tau _{m+i} \rightarrow A_{m+i}^{0} = 1, $
$ \vartheta _{l,p} \mapsto q_{l,p}^-, $
$i \in \{ 1, \ldots , n - 1 \} , $
$l \in \{ 1, \ldots , m \} , $
$l < p$
(where now we use the imbedding of \cite{Lawrence1990}).
By the braid valued Burau matrices and the similarity
of the recursion matrix with these,
the matrices
$ A_l^{(m)} $
for
$ l \in \{ 1+n, \ldots , n+m-1 \} $
represent
$B_n$,
if the matrices
$ A_l^{(m-1)}, $
$ l \in \{ 1+m, \ldots , n+m-1 \} $
and
$ b_{m,p}^{(m-1)}, $
$ p \in \{ 1+m, \ldots , n+m-1 \} $
represent
$B_nF_n$.
So finally we are led to a representation of
$B_n F_n \ldots F_{n+m-1}$
which is given by $\bar{\chi }$.

\bibliographystyle{amsalpha}

\end{document}